# A UNIFORM FUNCTIONAL LAW OF THE LOGARITHM FOR THE LOCAL EMPIRICAL PROCESS

By David M. Mason[1]

*University of Delaware*

We prove a uniform functional law of the logarithm for the local empirical process. To accomplish this we combine techniques from classical and abstract empirical process theory, Gaussian distributional approximation and probability on Banach spaces. The body of techniques we develop should prove useful to the study of the strong consistency of $d$-variate kernel-type nonparametric function estimators.

**1. Introduction.** Let $U, U_1, U_2, \ldots$, be a sequence of independent Uniform $[0, 1]$ random variables. Consider for each integer $n \geq 1$ the empirical distribution function based on $U_1, \ldots, U_n$,

$$G_n(t) = n^{-1} \sum_{i=1}^{n} 1\{U_i \leq t\}, \qquad -\infty < t < \infty.$$

Stute (1982a) was the first to initiate a concerted study of the almost sure behavior of the oscillation modulus of the uniform empirical process, which for any positive $0 < h < 1$ is defined to be

$$\varpi_n(h) = \sup\{|\sqrt{n}\{G_n(t+s) - G_n(t) - s\}| : 0 \leq t,\ t+s \leq 1,\ 0 \leq s \leq h\}.$$

He proved that whenever $\{h_n\}_{n \geq 1}$ is a sequence of positive constants converging to zero at a certain rate [see (H.i–iii)], then the following uniform law holds:

(1.1) $$\lim_{n \to \infty} \varpi_n(h_n)/\sqrt{2h_n \log(1/h_n)} = 1 \qquad \text{a.s.}$$

Now more generally, let $Z, Z_1, Z_2, \ldots$, be i.i.d. random variables taking values in $\mathbb{R}$ with common Lebesgue density function $f$. Stute (1982b) obtained

Received June 2001; revised April 2003.
[1]Supported in part by NSA Grant MDA904-02-1-0034 and NSF Grant DMS-02-03865.

from (1.1) and the probability integral transformation a uniform strong law for the kernel density estimator $f_n$ over compact intervals $J$; namely, he showed that under certain regularity conditions on $f$,

$$(1.2) \quad \lim_{n\to\infty} \sqrt{nh_n} \sup_{z\in J}\{|f_n(z) - Ef_n(z)|/\sqrt{2\|K\|_2^2 f(z)\log(1/h_n)}\} = 1 \quad \text{a.s.,}$$

where $f_n$ is defined, for $z \in \mathbb{R}$, to be

$$(1.3) \quad f_n(z) = (nh_n)^{-1} \sum_{i=1}^{n} K(h_n^{-1}(z - Z_i)),$$

with $K$ being a kernel with compact support and of bounded variation satisfying

$$(1.4) \quad 0 < \int_{\mathbb{R}} K^2(x)\,dx = \|K\|_2^2 < \infty.$$

Later, Stute (1984a) established a version of his strong law (1.1) for certain oscillations of the empirical process based upon $Z_1, Z_2, \ldots$, i.i.d. $d$-dimensional random vectors with common Lebesgue density function $f$. He used it to derive precise results on the uniform consistency of the $d$-variate kernel density estimator, which is defined as in (1.3), but with the $h_n$ inside $K$ replaced by $h_n^{1/d}$.

Deheuvels and Mason (1992) extended the Stute (1982a) strong law (1.1) to a uniform functional law of the logarithm (UFLL) for the cluster of random increment functions on $[0,1]$,

$$(1.5) \quad \{\xi_n(t,\cdot) : 0 \le t \le 1 - h_n\},$$

where, for each $0 \le t \le 1 - h_n$, $\xi_n(t,\cdot)$ is the function defined on $[0,1]$,

$$(1.6) \quad \xi_n(t,s) = \sqrt{n/h_n}\{G_n(t+h_ns) - G_n(t) - h_ns\}, \qquad 0 \le s \le 1,$$

and applied it to obtain exact rates of strong consistency for a number of nonparametric density estimators. (See Corollary 3 for a statement of this result.) Motivated partially by their work, Einmahl and Mason (2000) developed techniques from general empirical process theory and combined them with methods from Deheuvels and Mason (1992) to establish the precise rate of strong consistency over compact intervals for certain kernel-type nonparametric function estimators. Their results improved upon those of Härdle, Janssen and Serfling (1988), who had obtained only approximate rates. As a byproduct, they were able to obtain the Stute (1982b) result (1.2) through an approach based upon viewing $\{f_n(z) - Ef_n(z) : z \in J\}$ as an empirical process indexed by the class of functions $\{h_n^{-1}K(h_n^{-1}(z - \cdot)) : z \in J\}$. They also pointed out that the $d$-variate version of (1.2) could be derived in the same way. Giné and Guillou (2002) have recently done this and proved the somewhat unexpected result that whenever $K$ is continuous on $\mathbb{R}^d$ with support contained in $[-1/2, 1/2]^d$, and satisfies some additional assumptions (see Example F.1), the density $f$ is uniformly continuous on $\mathbb{R}^d$ and $\{h_n\}_{n \geq 1}$ satisfies conditions (H.i–iii), then

$$
(1.7) \quad \begin{aligned}
&\lim_{n \to \infty} \sup_{z \in \mathbb{R}^d} \sqrt{nh_n} |f_n(z) - Ef_n(z)| / \sqrt{2\|K\|_2^2 \log(1/h_n)} \\
&\quad = \sup_{z \in \mathbb{R}^d} \sqrt{f(z)} \qquad \text{a.s.}
\end{aligned}
$$

[Somewhat earlier, Deheuvels (2000) proved a dimension 1 version of this result.] We will derive a UFLL version of (1.7) as a corollary to our main result.

The proof of the Deheuvels and Mason (1992) UFLL for (1.5) was strongly based on the Komlós, Major and Tusnády (KMT) (1975) *Wiener process approximation* to partial sums of i.i.d. Poisson random variables, coupled with a functional large deviation result for the Wiener process. Such a precise and powerful strong approximation as given by KMT does not exist in the general empirical process setting.

Our goal in this paper is to show how one can meld the techniques from classical and abstract empirical process theory, Gaussian distributional approximation and probability on Banach spaces to prove a UFLL for a general indexed by class of functions version of (1.5) formed by a sequence of i.i.d. random variables $Z_1, Z_2, \ldots$, taking values in $\mathbb{R}^d$ with common Lebesgue density function $f$. The basic ingredients of our approach, along with their sources, are the following:

1. Poissonization [Einmahl (1987), Deheuvels and Mason (1992) and Giné, Mason and Zaitsev (2003)].
2. The Talagrand (1994) exponential inequality for the empirical process indexed by functions.
3. Tight bounds for the absolute moment for the supremum of the empirical process under a uniform covering number bound [Einmahl and Mason (2000) and Giné and Guillou (2001)].
4. Gaussian distributional approximation of multivariate sums [Zaitsev (1987a, b)].
5. Functional large deviation results for stochastic processes [Arcones (2003, 2004)].

We shall see that our approach is powerful enough to obtain the Deheuvels and Mason (1992) UFLL (without the use of KMT) as a corollary of our main result. The methods and results developed in this paper should be of potential application to the investigation of the strong consistency of a variety of multivariate nonparametric function estimators. To see how to apply the UFLL for the increment functions of the uniform empirical process to obtain exact rates of strong consistency for a number of univariate nonparametric density estimators, refer to Section 3 of Deheuvels and Mason (1992). Our main results are stated in Section 2, several examples are detailed in Section 3 and all the proofs are given in Section 4.

**2. Main results.** Let $Z, Z_1, Z_2, \ldots$, be i.i.d. $d$-dimensional random vectors with common Lebesgue density function $f$. Throughout this paper $\mathcal{G}$ will denote a class of measurable real valued functions defined on $\mathbb{R}^d$, which have support contained in $I^d := [-1/2, 1/2]^d$ and are bounded by some $\kappa > 0$. Let $|\cdot|_2$ denote the usual Euclidean norm on $\mathbb{R}^d$. Assume that the class $\mathcal{G}$ satisfies:

(G.i) $\lim_{|w|_2 \to 0} \sup_{g \in \mathcal{G}} \int_{\mathbb{R}^d} [g(x) - g(x+w)]^2 \, dx = 0$;

(G.ii) $\lim_{\lambda \to 1} \sup_{g \in \mathcal{G}} \int_{\mathbb{R}^d} [g(x) - g(\lambda x)]^2 \, dx = 0$.

In addition, let $\mathcal{F}$ denote the class of functions formed from $\mathcal{G}$ satisfying:

(F.i) for each $\lambda \geq 1$, $z \in \mathbb{R}^d$ and $g \in \mathcal{G}$, $g(z - \cdot \lambda) \in \mathcal{F}$.

To avoid using outer probability measures in all of our statements, we impose the measurability assumption:

(F.ii) $\mathcal{F}$ is a pointwise measurable class; that is, there exists a countable subclass $\mathcal{F}_0$ of $\mathcal{F}$ such that we can find, for any function $g \in \mathcal{F}$, a sequence of functions $\{g_m\}$ in $\mathcal{F}_0$ for which $g_m(z) \to g(z)$, $z \in \mathbb{R}^d$. [See Example 2.3.4 in van der Vaart and Wellner (1996).]

Finally we shall require the following entropy condition on the class $\mathcal{F}$. For $\varepsilon > 0$, let $N(\varepsilon, \mathcal{F}) = \sup_Q N(\kappa \varepsilon, \mathcal{F}, d_Q)$, where the supremum is taken over all probability measures $Q$ on $(\mathbb{R}^d, \mathcal{B})$, $d_Q$ is the $L_2(Q)$-metric, and, as usual, $N(\varepsilon, \mathcal{F}, d_Q)$ is the minimal number of balls $\{g : d_Q(g, g') < \varepsilon\}$ of $d_Q$-radius $\varepsilon$ needed to cover $\mathcal{F}$. Assume that $\mathcal{F}$ satisfies the following uniform entropy condition:

(F.iii) for some $C_0 > 0$ and $\nu_0 > 0$, $N(\varepsilon, \mathcal{F}) \leq C_0 \varepsilon^{-\nu_0}$, $0 < \varepsilon < 1$.

Let $\{h_n\}_{n \geq 1}$ be a sequence of positive constants less than 1 converging to zero. Choose any $z \in \mathbb{R}^d$. The local empirical process at $z$ indexed by $g \in \mathcal{G}$ is defined to be

$$(2.1) \quad E_n(z, g) := (nh_n)^{-1/2} \sum_{i=1}^n \{g(h_n^{-1/d}(z - Z_i)) - Eg(h_n^{-1/d}(z - Z))\}.$$

Einmahl and Mason (1997, 1998) obtained central limit theorems, strong approximations and functional laws of the iterated logarithm for the local empirical process at a fixed $z$. [Mason (1988) had treated a special case of this process, which he called the tail uniform empirical process.] They showed how to apply their results to obtain the exact rate of pointwise consistency for a number of well-known nonparametric kernel-type function estimators. The definition of the local empirical process given by Einmahl and Mason (1997, 1998) is a bit more general in that the $h_n^{1/d}$ is replaced by a sequence of bi-measurable functions. It extends an earlier notion introduced by Deheuvels and Mason (1994).

It is our aim to study the uniform limiting behavior of this process as $z$ moves over a compact set $J$. Towards this end we introduce the following normed versions of $E_n$: For any $z \in \mathbb{R}^d$ and $g \in \mathcal{G}$, set

$$D_n(z,g) := E_n(z,g)/\sqrt{2\log(1/h_n)} \tag{2.2}$$

and if $f(z) > 0$, set

$$L_n(z,g) := E_n(z,g)/\sqrt{2\log(1/h_n)f(z)}. \tag{2.3}$$

We will assume that the sequence $\{h_n\}_{n \geq 1}$ converges to zero at the following rate:

(H.i) $h_n \searrow 0$, $nh_n \nearrow \infty$;
(H.ii) $nh_n/\log(1/h_n) \to \infty$;
(H.iii) $\log(1/h_n)/\log\log n \to \infty$.

Consider the inner product defined for $g_1, g_2 \in \mathcal{G}$ by

$$(g_1, g_2) := \int_{I^d} g_1(u)g_2(u)\,du. \tag{2.4}$$

Let $G_2(I^d)$ be the Hilbert subspace of $L_2(I^d)$ spanned by $\mathcal{G}$. Now let $\mathcal{S}$ denote its reproducing kernel Hilbert space generated by the inner product $(\cdot,\cdot)$. Applying Theorem 4D of Parzen (1961), the space $\mathcal{S}$ can be represented as follows: Let $l_\infty(\mathcal{G})$, denote the class of bounded functions on $\mathcal{G}$. For any $\xi \in G_2(I^d)$, denote $\varphi_\xi \in l_\infty(\mathcal{G})$ by $\varphi_\xi(g) := (g,\xi)$, $g \in \mathcal{G}$. Each $\varphi_\xi$ is uniquely defined by $\xi$ in the sense that $\varphi_{\xi_1} = \varphi_{\xi_2}$ if and only if $\xi_1 = \xi_2$, in $L_2(I^d)$. The space $\mathcal{S} = \{\varphi_\xi : \xi \in G_2(I^d)\}$ has the inner product

$$\langle \varphi_{\xi_1}, \varphi_{\xi_2}\rangle := (\xi_1, \xi_2). \tag{2.5}$$

Let $\mathcal{S}_0$ denote the unit ball in $\mathcal{S}$ and, for any $\vartheta \in \mathcal{S}_0$ and $\varepsilon > 0$, set

$$B_\varepsilon(\vartheta) = \{\psi \in l_\infty(\mathcal{G}) : \|\psi - \vartheta\|_\mathcal{G} < \varepsilon\}, \tag{2.6}$$

where for any class of functions $\mathcal{C}$ and $\psi \in l_\infty(\mathcal{C})$, the class of bounded functions on $\mathcal{C}$,

$$\|\psi\|_\mathcal{C} = \sup_{g \in \mathcal{C}} |\psi(g)|. \tag{2.7}$$

Finally, write for any $\varepsilon > 0$,

$$\mathcal{S}_0^\varepsilon = \Big\{\psi \in l_\infty(\mathcal{G}) : \inf_{\vartheta \in \mathcal{S}_0} \|\psi - \vartheta\|_\mathcal{G} < \varepsilon\Big\}. \tag{2.8}$$

Throughout this paper $J$ will denote a compact subset of $\mathbb{R}^d$ with *nonempty interior*. For any $\gamma > 0$, we set

$$J_\gamma = \Big\{x : \inf_{z \in J} |x - z|_2 \leq \gamma\Big\}. \tag{2.9}$$

Our UFLL for the local empirical process is given in the following theorem.

THEOREM 1. *In addition to assumptions* (G.i–ii), (F.i–iii) *and* (H.i–iii), *assume that, for some $\gamma > 0$, the density $f$ is continuous and positive on $J_\gamma$. Then with probability* 1:

(a) *for all $\varepsilon > 0$, there exists an $n(\varepsilon)$ such that, for each $n \geq n(\varepsilon)$, $\{L_n(z, \cdot) : z \in J\} \subset \mathcal{S}_0^\varepsilon$;*

(b) *for any $\vartheta \in \mathcal{S}_0$ and $\varepsilon > 0$, there is an $n(\vartheta, \varepsilon)$ such that, for all $n \geq n(\vartheta, \varepsilon)$, there is a $z_n \in J$ such that $L_n(z_n, \cdot) \in B_\varepsilon(\vartheta)$.*

REMARK 1. It has long been recognized that the polynomial covering number assumption (F.iii) is the natural condition to impose upon the indexing class, when studying the local behavior of the empirical process. For instance, when Alexander (1987) made the first steps toward the investigation of the increments of the empirical process in a general indexed by a class of sets framework, he considered classes of index sets, which satisfy (F.iii). Nolan and Pollard (1987) and Nolan and Marron (1989) pointed out how the assumption (F.iii) on the class $\mathcal{F}$ arises naturally when investigating the large sample behavior of the kernel density estimator via empirical process indexed by a class of functions theory. (See Example F.1.) Later, Rio (1994) found that (F.iii) was the right assumption to impose on $\mathcal{F}$ when he derived his local invariance principle for the uniform ($[0,1]^d$) empirical process indexed by a class of functions, and applied it to kernel density estimation; as did Einmahl and Mason (1987, 2000, 2003) in their treatment of local empirical processes, Giné and Guillou (2002) in their derivation of rates of strong consistency for multivariate kernel density estimators and Deheuvels and Mason (2004) in their construction of universal confidence bands for regression functions. Classes of functions satisfying (F.iii) play a featured role in Devroye and Lugosi's (2000) derivation of bounds in the $L_1$ error for certain kinds of density estimators. This assumption also plays a critical role in the work of Giné, Koltchinskii and Wellner (2003) on ratio limit theorems for empirical processes.

REMARK 2. Condition (F.iii) was imposed to ensure that the moment bound (4.21) given in Fact 5 holds uniformly over all the classes of indexing functions considered in the proof of Theorem 1. One may surmise that (F.iii) could be replaced by a less restrictive entropy assumption. However, it is not clear whether this is the case. For a closely related discussion of this assumption, as it pertains to a local Gaussian process version of Theorem 1, see Mason (2003).

REMARK 3. It is routine to modify the proof of Theorem 1 to show that it remains true when (F.iii) is replaced by the bracketing condition:

(F.iii)′ for some $C_0 > 0$ and $\nu_0 > 0$, $N_{[\cdot]}(\varepsilon, \mathcal{F}, L_2(P)) \leq C_0 \varepsilon^{-\nu_0}$, $0 < \varepsilon < 1$.

Refer to page 270 of van der Vaart (1998) for the definition of $N_{[\,\cdot\,]}(\epsilon, \mathcal{F}, L_2(P))$. Essentially all that one has to do is to substitute the use of Fact 5 below by Lemma 19.34 of van der Vaart (1998).

Theorem 1 should be compared with the following functional law of the *iterated* logarithm at a fixed $z \in \mathbb{R}^d$ that can be inferred from Corollary 1.1 of Einmahl and Mason (1997), namely, that with probability 1 the sequence of processes indexed by $g \in \mathcal{G}$,

$$\{E_n(z,g)/\sqrt{2\log\log(1/h_n)f(z)},\ g \in \mathcal{G}\},$$

is relatively compact in $l_\infty(\mathcal{G})$ with set of limit points equal to $\mathcal{S}_0$. We see that to describe the behavior of $E_n(z,g)$ at a fixed point $z \in \mathbb{R}^d$, the $\sqrt{\log(1/h_n)}$ norming must be replaced by $\sqrt{\log\log(1/h_n)}$.

The following corollary provides a UFLL version of the Giné and Guillou (2002) result cited in (1.7).

COROLLARY 1. *In addition to assumptions* (G.i–ii), (F.i–iii) *and* (H.i–iii), *assume the density $f$ is uniformly continuous on $\mathbb{R}^d$. Then with probability 1:*

(a) *for all $\varepsilon > 0$, there exists an $n(\varepsilon)$ such that, for each $n \geq n(\varepsilon)$, $\{D_n(z,\cdot) : z \in \mathbb{R}^d\} \subset \{\tau_0 \mathcal{S}_0 : z \in \mathbb{R}^d\}^\varepsilon$, where $\tau_0 = \sup_{z \in \mathbb{R}^d} \sqrt{f(z)}$;*

(b) *for any $z \in \mathbb{R}^d$, $\vartheta \in \mathcal{S}_0$ and $\varepsilon > 0$, there is an $n(\vartheta, z, \varepsilon)$ such that, for all $n \geq n(\vartheta, z, \varepsilon)$, there is a $z_n \in \mathbb{R}^d$ such that $D_n(z_n, \cdot) \in \sqrt{f(z)} B_\varepsilon(\vartheta)$ and $|z_n - z|_2 < \varepsilon$.*

To see how Corollary 1 implies the Giné and Guillou (2002) result (1.7), let $\mathcal{G} = \{K\}$, where $K$ is continuous on $\mathbb{R}^d$ with support contained in $[-1/2, 1/2]^d$, and satisfies the conditions of Example F.1. Then assumptions (G.i–ii), (F.i–iii) and (H.i–iii) are satisfied. In this case $\mathcal{S}_0 = \{\varphi_\xi : \xi = uK/\|K\|_2 \text{ for some } |u| \leq 1\}$ and clearly $\sup\{|\varphi_\xi(K)| : \varphi_\xi \in \mathcal{S}_0\} = \|K\|_2$, from which (1.7) readily follows from parts (a) and (b) of Corollary 1.

Further examples are detailed in Section 3.

**3. Examples.** What classes of functions satisfy conditions (G.i–ii)? Using continuity of the shift and scale operators in $L_2(\mathbb{R}^d)$, it is trivial to see that (G.i–ii) hold for any class of functions $\mathcal{G}$ on $\mathbb{R}^d$ which is the convex hull of a finite number of $L_2(\mathbb{R}^d)$ functions. Here are some important classes that satisfy conditions (G.i–ii).

EXAMPLE G.1. Consider the class of functions $\mathcal{G}_c = \{\mathbb{1}_C : C \text{ is convex, closed and contained in } I^d\}$. Choose any $\mathbb{1}_C \in \mathcal{G}_c$, $0 < r < 1$ and $w \in \mathbb{R}^d$ satisfying $|w|_2 < r$; then

$$\int_{\mathbb{R}^d}[\mathbb{1}_C(x) - \mathbb{1}_C(x+w)]^2\,dx = |C\Delta(C-w)| \leq |C^r\Delta C| + |C^r\Delta(C-w)|$$
$$= 2\{|C^r| - |C|\},$$

where $\Delta$ denotes symmetric difference, $C^r = \{x : |x-y|_2 < r \text{ for some } y \in C\}$ and $|A|$ signifies the Lebesgue measure of a measurable set $A$. Now by the Steiner formula [see page 14 of Stoyan, Kendall and Mecke (1995)], we can conclude that, there exists a positive constant $c_d$ such that, for all $C$ convex, closed and contained in $I^d$, and $0 < r < 1$, we have $|C^r| \leq |C| + c_d r$. Thus $2\{|C^r| - |C|\} \leq 2c_d r$, which easily implies that the class satisfies condition (G.i). Condition (G.ii) is also readily verified.

EXAMPLE G.2. Let $\mathcal{G}$ be a bounded equicontinuous class of functions on $\mathbb{R}^d$ with support in $I^d$. From the inequality

$$\int_{\mathbb{R}^d}[g(x) - g(x+w)]^2\,dx = \int_{\mathbb{R}^d}[g(x)\mathbb{1}_{I^d}(x) - g(x+w)\mathbb{1}_{I^d}(x+w)]^2\,dx$$
$$\leq 2\int_{\mathbb{R}^d}[g(x)\{\mathbb{1}_{I^d}(x) - \mathbb{1}_{I^d}(x+w)\}]^2\,dx$$
$$+ 2\int_{\mathbb{R}^d}[g(x) - g(x+w)]^2\mathbb{1}_{I^d}(x+w)\,dx,$$

it is straightforward to show that condition (G.i) holds using the fact that the class of functions $\mathcal{G}$ is bounded and uniformly continuous in combination with Example G.1. Condition (G.ii) is checked in the same way.

Notice that the class $\mathcal{G}_+ = \{ag + b\mathbb{1}_C : g \in \mathcal{G}, \mathbb{1}_C \in \mathcal{G}_c \text{ and } |a| + |b| \leq D\}$, where $0 < D < \infty$ and $\mathcal{G}$ is any class of functions as in Example G.2 satisfies conditions (G.i–ii).

What about classes of functions $\mathcal{F}$ that satisfy all the conditions (G.i–ii) and (F.i–iii)?

EXAMPLE F.1. Set $\mathcal{G} = \{K\}$, where $K$ is continuous with support in $I^d$. Furthermore, whenever $d = 1$, assume $K$ is of bounded variation on $\mathbb{R}$, and whenever $d \geq 2$, that $K$ is of the form $K(x) = \Phi(x^T A x)$, for some $d \times d$ matrix $A$ and bounded continuous real valued function $\Phi$ of bounded variation on $\mathbb{R}$. Obviously (G.i–ii) hold for $\mathcal{G}$. The class $\mathcal{F}_K = \{K(z - \cdot\lambda) : \lambda \geq 1 \text{ and } z \in \mathbb{R}^d\}$ satisfies (F.i) by construction and (F.iii) by the results in Section 5 of Nolan and Pollard (1987). Also (F.ii) is readily verified using continuity of $K$.

EXAMPLE F.2. Let $\mathcal{G}_R = \{\mathbb{1}_R : R \in \mathcal{R}\}$, where $\mathcal{R} = $ class of closed rectangles contained in $I^d$, or $\mathcal{G}_E = \{\mathbb{1}_E : E \in \mathcal{E}\}$, where $\mathcal{E}$ is the class of closed ellipsoids contained in $I^d$. Clearly, since $\mathcal{G}_R$ and $\mathcal{G}_E$ are subsets of $\mathcal{G}_c$, conditions (Gi–ii) hold for both classes. Set

$$\text{(3.1)} \qquad \mathcal{F}_R = \{\mathbb{1}_R(z - \cdot\lambda) : \mathbb{1}_R \in \mathcal{G}_R, \ \lambda \geq 1 \text{ and } z \in \mathbb{R}^d\}$$

and

$$\text{(3.2)} \qquad \mathcal{F}_E = \{\mathbb{1}_E(z - \cdot\lambda) : \mathbb{1}_E \in \mathcal{G}_E, \ \lambda \geq 1 \text{ and } z \in \mathbb{R}^d\}.$$

It is well known that both the set of all closed rectangles and the set of all closed ellipsoids in $\mathbb{R}^d$ form Vapnik–Červonenkis (V.C.) classes; therefore, both $\mathcal{F}_R$ and $\mathcal{F}_E$ clearly satisfy (F.i) and (F.ii). [See van der Vaart and Wellner (1996) for the definition of a V.C. class, along with exercise 9 on page 151 of their book.] Finally, (F.iii) is readily verified for both $\mathcal{F}_R$ and $\mathcal{F}_E$.

Observe that the class of functions

$$\mathcal{F}_+ = \{ag_1 + bg_2 : g_1 \in \mathcal{F}_K, \ g_2 \in \mathcal{F}_R \text{ and } |a| + |b| \leq D\},$$

where $0 < D < \infty$, is easily shown to satisfy (G.i–ii) and (F.i–iii). This class should suffice for most applications.

The following corollary provides a UFLL version of Theorem 2.1 of Stute (1984).

COROLLARY 2. *Let $\{h_n\}_{n \geq 1}$ satisfy* (H.i–iii) *and let $\mathcal{G}_R$ and $\mathcal{F}_R$ be defined as above. Assume that, for some $\gamma > 0$, the density $f$ is continuous and positive on $[a_1 - \gamma, b_1 + \gamma] \times \cdots \times [a_d - \gamma, b_d + \gamma]$, where $-\infty < a_i < b_i < \infty$, $i = 1, \ldots, d$. Then with probability* 1:

(a) *for all $\varepsilon > 0$, there exists an $n(\varepsilon)$ such that, for each $n \geq n(\varepsilon)$, $\{L_n(z, \cdot) : z \in J\} \subset \mathcal{S}_0^\varepsilon$, where $J = [a_1, b_1] \times \cdots \times [a_d, b_d]$ and*

$$\text{(3.3)} \qquad \mathcal{S}_0 = \bigg\{\varphi : \varphi(\mathbb{1}_R) = \int_{I^d} \mathbb{1}_R(x)\xi(x)\,dx \text{ for } \mathbb{1}_R \in \mathcal{G}_R$$

$$\text{with } \xi \text{ satisfying } \int_{I^d} \xi^2(x)\,dx \leq 1\bigg\};$$

(b) *for any $\vartheta \in \mathcal{S}_0$ and $\varepsilon > 0$, there is an $n(\vartheta, \varepsilon)$ such that, for all $n \geq n(\vartheta, \varepsilon)$, there is a $z_n \in J$ such that $L_n(z_n, \cdot) \in B_\varepsilon(\vartheta)$.*

Notice that it is readily checked that, for $\mathcal{S}_0$ in (3.3), $\sup_{\varphi \in \mathcal{S}_0} \sup_{\mathbb{1}_R \in \mathcal{G}_R} |\varphi(\mathbb{1}_R)| = 1$, which on account of parts (a) and (b) of Corollary 2 implies that $\lim_{n \to \infty} \sup_{\mathbb{1}_R \in \mathcal{G}_R} |L_n(\mathbb{1}_R)| = 1$, a.s.

We end this section by showing how the UFLL for the increment functions of the uniform empirical process given in Theorem 3.1 of Deheuvels and Mason (1992) can be derived from Theorem 1. First note that the proof of Theorem 1 shows that it remains true when $I^d$ is replaced by any compact $d$-dimensional cube. In particular, in dimension 1, Theorem 1 holds with $I$ replaced by $[0, 1]$. Next, the classes of functions

$$\mathcal{G} = \{\mathbb{1}_{[0,t]} : t \in [0,1]\} \quad \text{and} \quad \mathcal{F} = \{\mathbb{1}_{[0,t]}(z - \cdot \lambda) : t \in [0,1],\ \lambda \geq 1 \text{ and } z \in \mathbb{R}\}$$

are readily shown to satisfy (G.i–ii) and (F.i–iii), respectively. Furthermore, in this setup,

(3.4)
$$\mathcal{S}_0 = \Big\{\varphi : \varphi(\mathbb{1}_{[0,t]}) = \int_0^t \xi(x)\,dx \text{ for } t \in [0,1]$$
$$\text{with } \xi \text{ satisfying } \int_0^1 \xi^2(x)\,dx \leq 1 \Big\}.$$

Recalling the notation in (1.5) and (1.6), set for each $n \geq 1$ and $t \in [0, 1 - h_n]$,

$$\mathcal{E}_n(t, \cdot) = \xi_n(t, \cdot)/\sqrt{2\log(1/h_n)}.$$

Assume that $\{h_n\}_{n\geq 1}$ satisfies (H.i–iii). Clearly, when the underlying distribution function is Uniform $[0,1]$, we can apply Theorem 1 to infer that, for any choice of $0 < \gamma < 1/2$ and for each $\varepsilon > 0$, there exists an $n(\varepsilon)$ such that, for any $n \geq n(\varepsilon)$, $\{\mathcal{E}_n(t, \cdot) : t \in [\gamma, 1 - \gamma]\} \subset \mathcal{S}_0^\varepsilon$. Combining this with (1.1), which implies that

(3.5)
$$\lim_{\gamma \searrow 0} \lim_{n \to \infty} \varpi_n(h_n \gamma)/\sqrt{2h_n \log(1/h_n)} = 0 \quad \text{a.s.,}$$

we obtain from Theorem 1 the following corollary, which is Theorem 3.1 of Deheuvels and Mason (1992). [Alternatively, in place of (1.1), we could have proved (3.5) by an argument based on Inequality 1.]

COROLLARY 3. *Let $\{h_n\}_{n\geq 1}$ satisfy (H.i–iii). Then with probability 1:*

(a) *for all $\varepsilon > 0$, there exists an $n(\varepsilon)$ such that, for each $n \geq n(\varepsilon)$, $\{\mathcal{E}_n(t, \cdot) : t \in [0, 1 - h_n]\} \subset \mathcal{S}_0^\varepsilon$ and*

(b) *for any $\vartheta \in \mathcal{S}_0$, $[a, b] \subset [0, 1]$, with $a < b$, and $\varepsilon > 0$, there is an $n(\vartheta, \varepsilon)$ such that, for all $n \geq n(\vartheta, \varepsilon)$, there is a $t_n \in [a, b]$ such that $\mathcal{E}_n(t_n, \cdot) \in B_\varepsilon(\vartheta)$.*

**4. Proofs.**

4.1. *Proof of part* (b) *of Theorem* 1.

4.1.1. *A large deviation result.* Set, for $n \geq 1$,

(4.1)
$$b_n = \sqrt{2nh_n \log(1/h_n)}.$$

Crucial to our proof is the following uniform large deviation result. Let $\eta_n$ be a standard Poisson random variable with rate $n$, independent of $Z, Z_1, Z_2, \ldots$, and consider the Poissonized version of the $L_n$ process:

$$\Pi_n(z,g) := (b_n\sqrt{f(z)})^{-1} \sum_{i=1}^{\eta_n} (g(h_n^{-1/d}(z-Z_i)) - Eg(h_n^{-1/d}(z-Z))),$$

where the empty sum is defined to be zero.

Define the rate function $I(\cdot)$ on $l_\infty(\mathcal{G})$ as follows. For any $\psi \in l_\infty(\mathcal{G})$,

(4.2) $$I(\psi) = \begin{cases} \frac{1}{2}\int_{I^d} \xi^2(u)\, du, & \text{if } \psi = \varphi_\xi \text{ for some } \xi \in G_2(I^d), \\ \infty, & \text{otherwise.} \end{cases}$$

Recall the definitions of $\varphi_\xi$ and $G_2(I^d)$ between (2.4) and (2.5) in Section 2. Also denote for any subset $B \subset l_\infty(\mathcal{G})$,

(4.3) $$I(B) = \inf\{I(\psi) : \psi \in B\}.$$

We endow $l_\infty(\mathcal{G})$ with the topology generated by the norm $\|\cdot\|_\mathcal{G}$, defined as in (2.7).

PROPOSITION 1. *Under the assumptions of Theorem* 1, *for any sequence* $\{m_n\}_{n\geq 1}$ *of positive integers and any triangular array of points* $z_{i,n}$, $i = 1, \ldots, m_n$, $n \geq 1$, *in $J$, we have:*

(i) *for all closed subsets $F$ of $l_\infty(\mathcal{G})$,*

$$\limsup_{n\to\infty} \max_{1\leq i\leq m_n} \varepsilon_n \log P\{\Pi_{i,n}(\cdot) \in F\} \leq -I(F);$$

(ii) *for all open subsets $G$ of $l_\infty(\mathcal{G})$,*

$$\liminf_{n\to\infty} \min_{1\leq i\leq m_n} \varepsilon_n \log P\{\Pi_{i,n}(\cdot) \in G\} \geq -I(G),$$

*where $\Pi_{i,n}(\cdot) = \Pi_n(z_{i,n}, \cdot)$, $i = 1, \ldots, m_n$, $n \geq 1$ and*

(4.4) $$\varepsilon_n = (2\log(1/h_n))^{-1}.$$

4.1.2. *Proof of Proposition* 1. We will take advantage of some recent work of Arcones (2003, 2004). In fact, we shall require the following trivial generalization of Theorem 3.1 of Arcones (2003). In the statement of this result, $P^*$ and $P_*$ denote the usual outer and inner measures associated with $P$, and $\overline{A}$ and $A^o$ denote the closure and interior of $A$, respectively. Let $l_\infty(T)$ denote the space of bounded functions on $T$. Also LDP is short for *large deviation principle*, as defined, for instance, in Arcones (2003). Note that a basic ingredient of Fact 1 is the uniform exponential tightness condition (A.iii).

FACT 1. *Let $\{X_{i,n}(t) : t \in T, 1 \leq i \leq m_n\}$, where $\{m_n\}_{n \geq 1}$ is a sequence of positive integers, be a triangular array of stochastic processes and let $T$ be an index set. Let $\{\varepsilon_n\}_{n \geq 1}$ be a sequence of positive numbers that converges to zero. Let $\varrho$ be a pseudo-metric on $T$. Consider the following conditions:*

(A.i) *$(T, \varrho)$ is totally bounded.*

(A.ii) *For each choice of $t_1, \ldots, t_k \in T$, the triangular array of vectors $\{(X_{i,n}(t_1), \ldots, X_{i,n}(t_k)), 1 \leq i \leq m_n\}$ satisfies uniformly the LDP with speed $\varepsilon_n$ and good rate function $I_{t_1,\ldots,t_k}$, in the sense that, for any Borel subset $A \subset \mathbb{R}^k$,*

$$-\inf_{z \in A^o} I_{t_1,\ldots,t_k}(z) \leq \liminf_{n \to \infty} \varepsilon_n \min_{1 \leq i \leq m_n} \log P_*\{(X_{i,n}(t_1), \ldots, X_{i,n}(t_k)) \in A\}$$

$$\leq \limsup_{n \to \infty} \varepsilon_n \max_{1 \leq i \leq m_n} \log P^*\{(X_{i,n}(t_1), \ldots, X_{i,n}(t_k)) \in A\}$$

$$\leq -\inf_{z \in \overline{A}} I_{t_1,\ldots,t_k}(z),$$

*and for any $0 < \alpha < \infty$, the set $\{z \in \mathbb{R}^k : I_{t_1,\ldots,t_k}(z) \leq \alpha\}$ is a compact set in $\mathbb{R}^k$.*

(A.iii) *For each $\tau > 0$,*

$$\lim_{\eta \to 0} \limsup_{n \to \infty} \varepsilon_n \max_{1 \leq i \leq m_n} \log P^*\left\{\sup_{\varrho(s,t) \leq \eta} |X_{i,n}(t) - X_{i,n}(s)| \geq \tau\right\} = -\infty.$$

*Then, for each $0 < \alpha < \infty$, the set $\{\psi \in l_\infty(T) : I(\psi) \leq \alpha\}$ is a compact set in $l_\infty(T)$, where*

$$I(\psi) = \sup\{I_{t_1,\ldots,t_k}(\psi(t_1), \ldots, \psi(t_k)) : t_1, \ldots, t_k \in T, \ k \geq 1\}.$$

*Moreover, one gets the following upper and lower bounds in the LDP with respect to outer and inner probabilities (because of possible lack of measurability): For each $A \subset l_\infty(T)$,*

$$-\inf_{\psi \in A^o} I(\psi) \leq \liminf_{n \to \infty} \varepsilon_n \min_{1 \leq i \leq m_n} \log P_*\{X_{i,n} \in A\}$$

$$\leq \limsup_{n \to \infty} \varepsilon_n \max_{1 \leq i \leq m_n} \log P^*\{X_{i,n} \in A\} \leq -\inf_{\psi \in \overline{A}} I(\psi).$$

Also we will require the following fact, which follows by applying Theorem 5.2 of Arcones (2004) to a finite index set $T$.

FACT 2. *Let $\{U_{i,n}(t) : t \in T, 1 \leq i \leq m_n\}$ be a triangular array of centered Gaussian random vectors, where $\{m_n\}_{n \geq 1}$ is a sequence of positive integers and $T = \{t_1, \ldots, t_d\}$ is a finite index set. Let $\{\varepsilon_n\}_{n \geq 1}$ be a sequence of positive numbers that converges to zero as $n \to \infty$. Assume that for a covariance matrix $R = \{R(t_i, t_j) : (t_i, t_j) \in T^2\}$, we have, for any $s, t \in T$,*

(4.5) $$\lim_{n \to \infty} \max_{1 \leq i \leq m_n} |R(s,t) - \varepsilon_n^{-1} E[U_{i,n}(s) U_{i,n}(t)]| = 0.$$

*Then for any Borel subset $A \subset \mathbb{R}^d$,*

$$
\begin{aligned}
-\inf_{z \in A^o} I_{t_1,\ldots,t_d}(z) &\leq \liminf_{n \to \infty} \varepsilon_n \min_{1 \leq i \leq m_n} \log P\{U_{i,n} \in A\} \\
&\leq \limsup_{n \to \infty} \varepsilon_n \max_{1 \leq i \leq m_n} \log P\{U_{i,n} \in A\} \leq -\inf_{z \in \overline{A}} I_{t_1,\ldots,t_d}(z),
\end{aligned}
\tag{4.6}
$$

*where for $z \in \mathbb{R}^d$,*

$$I_{t_1,\ldots,t_d}(z) = \inf\{2^{-1}\xi^T R\xi : R\xi = z\}. \tag{4.7}$$

The following lemma, which is a special case of a result of Stein [(1970), pages 62 and 63], will come in handy.

LEMMA 1. *Let $f$ be a Lebesgue density function on $\mathbb{R}^d$, which for some $\gamma > 0$ is bounded and uniformly continuous on $D_\gamma$, where $D$ is a closed subset of $\mathbb{R}^d$ and $D_\gamma$ is defined as in* (2.9). *Then for any $L_1(\mathbb{R}^d)$ function $H$, which is equal to zero for $x \notin I^d$,*

$$\sup_{z \in D} |f * H_h(z) - I(H)f(z)| \to 0 \qquad \text{as } h \searrow 0, \tag{4.8}$$

*where $I(H) = \int_{\mathbb{R}^d} H(u)\,du$ and $f * H_h(z) := h^{-1} \int_{\mathbb{R}^d} f(x) H(h^{-1/d}(z-x))\,dx$.*

Choose $\mathcal{G}_q = \{g_1, \ldots, g_q\} \subset \mathcal{G}$ and $z_{1,n}, \ldots, z_{m_n,n} \in J$. Let $\{U_{i,n}(g) : g \in \mathcal{G}_q, 1 \leq i \leq m_n\}$ be a triangular array of centered Gaussian random vectors each with covariance function

$$\sigma_{i,n}(g_l, g_k) = n(b_n^2 f(z_{i,n}))^{-1} \operatorname{cov}(g_l(h_n^{-1/d}(z_{i,n} - Z)), g_k(h_n^{-1/d}(z_{i,n} - Z))),$$

$$1 \leq i \leq m_n.$$

It is routine using Lemma 1 to show that, with $\varepsilon_n$ as in (4.4), as $n \to \infty$,

$$\max_{1 \leq i \leq m_n} \max_{1 \leq l,k \leq q} |\varepsilon_n^{-1} \sigma_{i,n}(g_l, g_k) - \sigma(g_l, g_k)| \to 0, \tag{4.9}$$

where $\sigma(g_l, g_k) := \int_{I^d} g_l(u) g_k(u)\,du$.

Thus Fact 2 applies here and its conclusion (4.6) holds with

$$R(g_l, g_k) = \sigma(g_l, g_k). \tag{4.10}$$

Consider now the triangular array of Poisson-type processes indexed by $g \in \mathcal{G}$,

$$
\begin{aligned}
\Pi_{i,n}(g) &:= (b_n \sqrt{f(z_{i,n})})^{-1} \\
&\quad \times \sum_{j=1}^{\eta_n} (g(h_n^{-1/d}(z_{i,n} - Z_j)) - Eg(h_n^{-1/d}(z_{i,n} - Z))),
\end{aligned}
\tag{4.11}
$$

$1 \le i \le m_n$, $n \ge 1$. Notice that for each $1 \le i \le m_n$, the process $\{\Pi_{i,n}(g)\}_{g \in \mathcal{G}_q}$ has the same covariance function as the process $\{U_{i,n}(g)\}_{g \in \mathcal{G}_q}$. We claim that for any Borel subset $A \subset \mathbb{R}^q$,

$$
\begin{aligned}
-\inf_{z \in A^o} I_{g_1,\ldots,g_q}(z) &\le \liminf_{n \to \infty} \varepsilon_n \min_{1 \le i \le m_n} \log P\{\Pi_{i,n} \in A\} \\
&\le \limsup_{n \to \infty} \varepsilon_n \max_{1 \le i \le m_n} \log P\{\Pi_{i,n} \in A\} \\
&\le -\inf_{z \in \overline{A}} I_{g_1,\ldots,g_q}(z),
\end{aligned}
\tag{4.12}
$$

where $\varepsilon_n$ is as in (4.4), and for $z \in \mathbb{R}^q$, $I_{g_1,\ldots,g_q}(z)$ is defined as in (4.7), with $R$ as in (4.10).

To show this we shall need the following result of Zaitsev (1987a). For probability measures $P$ and $Q$ on the Borel subsets of $\mathbb{R}^q$, $q \ge 1$, and $\delta > 0$, let

$$
(4.13) \quad \lambda(P, Q, \delta) := \sup\{P(A) - Q(A^\delta),\ Q(A) - P(A^\delta) : A \subset \mathbb{R}^q,\ \text{Borel}\},
$$

where $A^\delta$ denotes the $\delta$-neighborhood of $A$, $A^\delta := \{x \in \mathbb{R}^q : \inf_{y \in A} |x - y|_2 < \delta\}$. Let $P$ be an infinitely divisible $q$-dimensional distribution with spectral measure concentrated on the ball $\{x \in \mathbb{R}^q : |x|_2 \le \beta\}$, $\beta > 0$, and let $Q$ be the $q$-dimensional normal distribution with the same mean and covariance matrix as $P$. The following inequality is contained in Theorem 1.1 and Example 1.2 of Zaitsev (1987a). See, as well, a slightly weaker statement in Theorem 1.2 of Zaitsev (1987b).

FACT 3. *For all $\delta \ge 0$,*

$$
(4.14) \quad \lambda(P, Q, \delta) \le c_{1,q} \exp(-\delta/(c_{2,q}\beta)),
$$

*where $c_{i,q} \le c_i q^2$ with $c_1$, $c_2$ being universal finite positive constants.*

It is easy to see that the distribution of $(\Pi_{i,n}(g_1), \ldots, \Pi_{i,n}(g_q))$, being compound Poisson, is infinitely divisible with spectral measure, uniformly in $1 \le i \le m_n$, concentrated on the ball $\{x \in \mathbb{R}^q : |x|_2 \le \beta\}$, where $\beta = \rho/\sqrt{nh_n \log(1/h_n)}$ and $\rho > 0$ is a constant. This follows from the fact that, for some $\rho' > 0$, uniformly in $1 \le i \le m_n$, $n \ge 1$,

$$
\begin{aligned}
(b_n\sqrt{f(z_{i,n})})^{-1}&|g(h_n^{-1/d}(z_{i,n} - Z_j)) - Eg(h_n^{-1/d}(z_{i,n} - Z))| \\
&\le \rho'/\sqrt{nh_n \log(1/h_n)}.
\end{aligned}
$$

Hence by applying Fact 3, we see that uniformly in $1 \le i \le m_n$,

$$
(4.15) \quad \lambda(\mathbf{P}_{i,n}, \mathbf{Q}_{i,n}, \delta) \le c_{1,q} \exp(-\delta\sqrt{nh_n \log(1/h_n)}/(\rho c_{2,q})),
$$

where $\mathbf{P}_{i,n}$ is the distribution of $(\Pi_{i,n}(g_1), \ldots, \Pi_{i,n}(g_q))$ and $\mathbf{Q}_{i,n}$ is the distribution of $(U_{i,n}(g_1), \ldots, U_{i,n}(g_q))$. Using (H.ii), it is easy to infer from (4.15) that, for each $\delta > 0$, as $n \to \infty$,

$$(4.16) \qquad \max_{1 \leq i \leq m_n} \log(1/h_n)^{-1} \log \lambda(\mathbf{P}_{i,n}, \mathbf{Q}_{i,n}, \delta) \to -\infty.$$

Therefore, since by (4.9) in combination with Fact 2, (4.12) holds with $\Pi_{i,n}$ replaced by $U_{i,n}$, we readily conclude from (4.16) that (4.12) is satisfied. Hence assumption (A.ii) of Fact 1 holds with $\varepsilon_n$ as in (4.4).

Our next goal is to verify that (A.iii) holds with $T = \mathcal{G}$. Let $\varrho$ denote the pseudo-metric on $\mathcal{G}$,

$$(4.17) \qquad \varrho(g, g') = \sqrt{\int_{\mathbb{R}^d} (g(u) - g'(u))^2 \, du}, \qquad g, g' \in \mathcal{G}.$$

We shall show that for each $\tau > 0$,

$$\lim_{\eta \to 0} \limsup_{n \to \infty} \varepsilon_n \max_{1 \leq i \leq m_n} \log P^* \left\{ \sup_{\varrho(g,g') \leq \eta} |\Pi_{i,n}(g) - \Pi_{i,n}(g')| \geq \tau \right\} = -\infty.$$

Observe that whenever $\eta_n = m$, for some $m \geq 1$, for $g, g' \in \mathcal{G}$ and each $1 \leq i \leq m_n$,

$$\Pi_{i,n}(g) - \Pi_{i,n}(g')$$
$$= (b_n \sqrt{f(z_{i,n})})^{-1} \sum_{j=1}^{m} \{(g - g')(h_n^{-1/d}(z_{i,n} - Z_j))$$
$$\qquad - E(g - g')(h_n^{-1/d}(z_{i,n} - Z))\}$$
$$=: (b_n \sqrt{f(z_{i,n})})^{-1} T_{m,n,i}(g - g').$$

Notice that we can choose $0 < \beta_1 < \beta_2 < \infty$ such that, for all $z \in J_\gamma$,

$$(4.18) \qquad 0 < \beta_1 \leq f(z) \leq \beta_2 < \infty.$$

Thus for each $1 \leq i \leq m_n$,

$$P\left\{ \sup_{\varrho(g,g') \leq \eta} |\Pi_{i,n}(g) - \Pi_{i,n}(g')| \geq \tau/\sqrt{\beta_1} \right\}$$
$$(4.19) \qquad \leq P\left\{ \max_{1 \leq i \leq m_n} \max_{1 \leq m \leq 2n} \sup_{\varrho(g,g') \leq \eta} |T_{m,n,i}(g - g')| \geq \tau b_n \right\} + P\{\eta_n > 2n\}$$
$$:= p_n(\tau, \eta) + P\{\eta_n > 2n\}.$$

To finish the proof, we shall require two more facts and an inequality following from them.

Let $X, X_1, X_2, \ldots$, be i.i.d. on a probability space $(\mathcal{X}, \mathcal{A}, P)$. Let $\mathcal{H}$ be a pointwise measurable class of real valued functions defined on $\mathcal{X}$. Further let $\epsilon_1, \epsilon_2, \ldots$, be a sequence of i.i.d. Rademacher random variables independent of $X_1, X_2, \ldots$. [By Rademacher, we mean $P(\epsilon_1 = 1) = P(\epsilon_1 = -1) = 1/2$.] Set, for each $g \in \mathcal{H}$ and $m \geq 1$,

$$T_m(g) = \sum_{j=1}^{m} \{g(X_j) - Eg(X)\}. \tag{4.20}$$

We shall need the following inequality, which is essentially due to Talagrand (1994). See Einmahl and Mason (2000).

FACT 4. *Let $\mathcal{H}$ be a pointwise measurable class of functions on $(\mathcal{X}, \mathcal{A})$ satisfying, for some $0 < M < \infty$, $\|g\|_{\mathcal{X}} := \sup_{x \in \mathcal{X}} |g(x)| \leq M$, $g \in \mathcal{H}$. Then, for all $n \geq 1$ and $t > 0$, we have, for suitable finite constants $A_1, A_2 > 0$,*

$$P\left\{\max_{1 \leq m \leq n} \|T_m\|_{\mathcal{H}} \geq A_1\left(E\left\|\sum_{i=1}^{n} \epsilon_i g(X_i)\right\|_{\mathcal{H}} + t\right)\right\}$$
$$\leq 2[\exp(-A_2 t^2 / n\sigma_{\mathcal{H}}^2) + \exp(-A_2 t/M)],$$

*where $\sigma_{\mathcal{H}}^2 = \sup_{g \in \mathcal{H}} \mathrm{Var}(g(X))$.*

Let $G$ be a finite valued measurable function satisfying, for all $x \in \mathcal{X}$, $G(x) \geq \sup_{g \in \mathcal{H}} |g(x)|$, and define

$$N(\varepsilon, \mathcal{H}) = \sup_{Q} N(\varepsilon \sqrt{Q(G^2)}, \mathcal{H}, d_Q),$$

where the supremum is taken over all probability measures $Q$ on $(\mathcal{X}, \mathcal{A})$ for which $0 < Q(G^2) < \infty$ and $d_Q$ is the $L_2(Q)$-metric. As above, $N(\varepsilon, \mathcal{H}, d_Q)$ is the minimal number of balls $\{g : d_Q(g, g') < \varepsilon\}$ of $d_Q$-radius $\varepsilon$ needed to cover $\mathcal{H}$.

We shall require the following moment bound of Einmahl and Mason (2000). [For a similar bound, refer to Giné and Guillou (2001).]

FACT 5. *Let $\mathcal{H}$ be a pointwise measurable class of real valued bounded functions on $(\mathcal{X}, \mathcal{A})$ such that, for some constants $\beta, \nu, C > 1$, $\sigma \leq 1/(8C)$ and function $G$ as above, the following four conditions hold:*

$$Q(G^2) = EG^2(X) \leq \beta^2;$$
$$N(\varepsilon, \mathcal{H}) \leq C\varepsilon^{-\nu}, \qquad 0 < \varepsilon < 1;$$
$$\sigma_0^2 := \sup_{g \in \mathcal{H}} Eg^2(X) \leq \sigma^2;$$

*and*

$$\sup_{g \in \mathcal{H}} \|g\|_{\mathcal{X}} \leq (2\sqrt{\nu+1})^{-1}\sqrt{n\sigma^2/\log(\beta \vee 1/\sigma)}.$$

*Then we have, for a universal constant $A_3$,*

(4.21) $$E\left\|\sum_{i=1}^{n} \epsilon_i g(X_i)\right\|_{\mathcal{H}} \leq A_3\sqrt{\nu n \sigma^2 \log(\beta \vee 1/\sigma)}.$$

We shall make frequent use of the next inequality, which follows readily from Facts 4 and 5.

INEQUALITY 1. *Let $\{\mathcal{F}_n\}_{n \geq 1}$ be a sequence of classes of measurable real valued functions on $\mathbb{R}^d$ each bounded by $M > 0$ and satisfying uniformly in $n \geq 1$, (F.ii–iii). Let $\{h_n\}_{n \geq 1}$ be a sequence of positive constants less than 1 converging to zero at the rate (H.ii). Assume that, for some $\gamma > 0$, for all $n$ large,*

(4.22) $$\sigma_n^2 := \sup_{g \in \mathcal{F}_n} Eg^2(X) \leq \gamma^2 h_n.$$

*Then with $b_n$ as in (4.1), there exist constants $D_0 > 0$ and $D_1 > 0$ such that, for all $\rho > 0$ and all $n$ large,*

$$P_n(\rho) := P\left\{\max_{1 \leq m \leq 2n} \|T_m\|_{\mathcal{F}_n} \geq (\gamma + \rho)D_1 b_n\right\}$$

(4.23) $$\leq 2\exp(-D_0(\rho/\gamma)^2 \log(1/h_n)).$$

PROOF. First, by Fact 4, for suitable finite constants $A_1, A_2 > 0$,

$$P\left\{\max_{1 \leq m \leq 2n} \|T_m\|_{\mathcal{F}_n} \geq A_1\left(E\left\|\sum_{i=1}^{2n}\epsilon_i g(X_i)\right\|_{\mathcal{F}_n} + \rho b_n\right)\right\}$$
$$\leq 2[\exp(-A_2\rho^2 b_n^2/(2n\sigma_n^2)) + \exp(-A_2\rho b_n/M)].$$

Now by using (4.22) with Fact 5, we get, for a suitable $D_1 \geq A_1$, for all $n$ large,

$$A_1 E\left\|\sum_{i=1}^{2n}\epsilon_i g(X_i)\right\|_{\mathcal{F}_n} \leq \gamma D_1 b_n,$$

which gives

$$P_n(\rho) \leq P\left\{\max_{1 \leq m \leq 2n} \|T_m\|_{\mathcal{F}_n} \geq A_1\left(E\left\|\sum_{i=1}^{2n}\epsilon_i g(X_i)\right\|_{\mathcal{F}_n} + \rho b_n\right)\right\}.$$

Therefore we readily conclude from Fact 4 that, for some constants $D_2, D_3 > 0$,

$$P_n(\rho) \leq 2\exp(-D_2(\rho/\gamma)^2 \log(1/h_n)) + 2\exp(-D_3\rho b_n),$$

which by (H.ii) is, for some $D_0 > 0$ and all large $n$, less than or equal to

$$2\exp(-D_0(\rho/\gamma)^2 \log(1/h_n)). \qquad \square$$

Returning to the proof of Proposition 1, for any $\eta > 0$, with $\varrho$ as in (4.17), let

$$(4.24) \qquad \mathcal{H}(\eta) := \{g - g' : \varrho(g, g') \leq \eta,\ g, g' \in \mathcal{G}\}$$

and

$$(4.25) \quad \mathcal{H}_n(\eta) = \{(g - g')(h_n^{-1/d}(z_{i,n} - \cdot)) : g - g' \in \mathcal{H}(\eta), 1 \leq i \leq m_n\}.$$

Using this notation, we can write $p_n(\tau, \eta)$ in (4.19) as

$$(4.26) \qquad p_n(\tau, \eta) = P\left\{\max_{1 \leq m \leq 2n} \|T_m\|_{\mathcal{H}_n(\eta)} \geq \tau b_n\right\}.$$

Clearly, by using the fact that each $g - g'$ is uniformly bounded by $M := 2\kappa$, we get, with $D_1$ as in Inequality 1 that for some $\delta > 0$ and all $\eta > 0$,

$$(4.27) \qquad \max_{1 \leq i \leq m_n} \sup_{g - g' \in \mathcal{H}(\eta)} E(g - g')^2(h_n^{-1/d}(z_{i,n} - Z)) \leq h_n \delta^2 \eta^2 / D_1^2.$$

Let $\mathcal{F}' = \{g - g' : g, g' \in \mathcal{F}\}$. Clearly, from (F.iii) we get that, for some $C > 0$ and with $\nu = 2\nu_0$, $N(\varepsilon, \mathcal{F}') \leq C\varepsilon^{-\nu}$, $0 < \varepsilon < 1$.

Now for any $\tau > 0$ and $\eta > 0$ such that $\eta\delta < \tau/2$, we see that

$$p_n(\tau, \eta) \leq P\left\{\max_{1 \leq m \leq 2n} \|T_m\|_{\mathcal{H}_n(\eta)} \geq (\eta\delta + \tau/2)b_n\right\}.$$

Therefore, by (4.27) and $\mathcal{H}_n(\eta) \subset \mathcal{F}'$, we can apply Inequality 1 with $\gamma = \eta\delta/D_1$ and $\rho = \tau/(2D_1)$ to show that, for all $\tau > 0$ and sufficiently small $\eta > 0$ satisfying $\eta\delta < \tau/2$ and all $n$ large enough,

$$(4.28) \qquad p_n(\tau, \eta) \leq 2\exp(-4^{-1}\delta^{-2}D_0(\tau/\eta)^2 \log(1/h_n)),$$

which implies that

$$\limsup_{n \to \infty} \log p_n(\tau, \eta)/(2\log(1/h_n)) \leq -8^{-1}\delta^{-2}D_0(\tau/\eta)^2.$$

Next, by Chebyshev's inequality applied to $\exp(\eta_n \log 2)$, we get $P\{\eta_n > 2n\} \leq \exp((1 - 2\log 2)n)$, from which we obtain that

$$\lim_{n \to \infty} \log P\{\eta_n > 2n\}/(2\log(1/h_n)) = -\infty.$$

Putting everything together, taking (4.19) into account, we conclude with $\varepsilon_n = (2\log(1/h_n)t)^{-1}$,

$$\lim_{\eta \to 0} \limsup_{n \to \infty} \varepsilon_n \max_{1 \leq i \leq m_n} P\bigg\{\sup_{\varrho(g,g') \leq \eta} |\Pi_{i,n}(g) - \Pi_{i,n}(g')| \geq \tau/\sqrt{\beta_1}\bigg\}$$
$$\leq \lim_{\eta \to 0}(-8^{-1}\delta^{-2}D_0(\tau/\eta)^2) = -\infty.$$

This shows that condition (A.iii) of Fact 1 holds. Assumption (F.iii) implies that (A.i) of Fact 1 is satisfied and we have already verified (A.ii) above. Thus we have checked all the conditions of Fact 1 and can infer that its conclusions hold for the triangular array of processes $\{\Pi_{i,n}(\cdot), i=1,\ldots,m_n\}$. Finally, it can be deduced from Theorem 4.2 of Arcones (2004) that, in our situation, the $I(\psi)$ as arises from Fact 1 has the representation (4.2). This completes the proof of Proposition 1.

4.1.3. *Poissonization.* Choose $z_{1,n},\ldots,z_{m_n,n} \in J$, and, for $1 \leq i \leq m_n$, let, for $g \in \mathcal{G}$,

$$L_{i,n}(g) := (b_n\sqrt{f(z_{i,n})})^{-1} \sum_{j=1}^n (g(h_n^{-1/d}(z_{i,n} - Z_j)) - Eg(h_n^{-1/d}(z_{i,n} - Z))).$$

We shall need the following special case of Lemma 2.1 of Giné, Mason and Zaitsev (2003). Its idea may be traced back to Pyke and Shorack [(1968), proof of Lemma 2.2], through Einmahl (1987) and Deheuvels and Mason (1992) [also see Einmahl and Mason (1997)]. For a further generalization, along with additional historical remarks, refer to Borisov (2002).

FACT 6. *Choose $z_{1,n},\ldots,z_{m_n,n} \in J$. Whenever*

(4.29) $$\sum_{i=1}^{m_n} P\{Z \in z_{i,n} - h_n^{1/d}I^d\} \leq 1/2,$$

*then for all Borel subsets $B_1,\ldots,B_{m_n}$ of $l_\infty(\mathcal{G})$,*

$$P\{L_{i,n} \in B_i, i=1,\ldots,m_n\} \leq 2P\{\Pi_{i,n} \in B_i, i=1,\ldots,m_n\},$$

*where $\Pi_{i,n}$ is the Poissonized version of $L_{i,n}$ as defined in* (4.11).

Our next lemma completes the proof of part (b) of Theorem 1.

LEMMA 2. *With probability 1, for any $\vartheta \in \mathcal{S}_0$, and $\varepsilon > 0$, there is an $n(\vartheta,\varepsilon)$ such that, for all $n \geq n(\vartheta,\varepsilon)$, there is a $z_n \in J$ such that*

(4.30) $$L_n(z_n,\cdot) \in B_\varepsilon(\vartheta).$$

PROOF. Recall definitions (4.2) and (4.3). Choose any $\vartheta \in \mathcal{S}_0$ with $0 < \langle \vartheta, \vartheta \rangle = 2I(\vartheta) \leq 1$ and $\varepsilon > 0$ small enough so that

(4.31) $$0 < 2I(B_\varepsilon(\vartheta)) < 1.$$

Select $z_{1,n}, \ldots, z_{m_n,n} \in J$ such that the components of $z_{i,n}$ and $z_{j,n}$, $i \neq j$, differ in absolute value by more than $h_n^{1/d}$, (4.29) holds and

(4.32) $$\log m_n / \log(1/h_n) \to 1 \quad \text{as } n \to \infty.$$

The existence of such a sequence $\{m_n\}_{n \geq 1}$ is guaranteed by (4.18) and the assumption that $J$ has nonempty interior, which implies that $[a_1, b_1] \times \cdots \times [a_d, b_d] \subset J$, for some $-\infty < a_i < b_i < \infty$, $i = 1, \ldots, d$.

We see by Fact 6 that

$$P_n = P\{L_{i,n} \notin B_\varepsilon(\vartheta),\ i = 1, \ldots, m_n\} \leq 2P\{\Pi_{i,n} \notin B_\varepsilon(\vartheta),\ i = 1, \ldots, m_n\}.$$

Now by using the independence property of the Poisson processes $\Pi_{i,n}$, $1 \leq i \leq m_n$, following from the choice of the $z_{i,n}$, $i = 1, \ldots, m_n$, and the assumption that the functions $g$ have support in $I^d$, this last bound equals

$$2 \prod_{i=1}^{m_n} P\{\Pi_{i,n} \notin B_\varepsilon(\vartheta)\}.$$

Applying part (ii) of Proposition 1, we see that this last expression is, for any $\rho > 0$ and all $n$ sufficiently large,

$$\leq 2[1 - \exp(-2(1 + \rho) \log(1/h_n) I(B_\varepsilon(\vartheta)))]^{m_n},$$

which, in turn, by (4.31) and an appropriate choice of $0 < \rho < 1$ is, for some $0 < \tau < 1$ and for all $n$ sufficiently large,

$$\leq 2[1 - \exp(-\tau \log(1/h_n))]^{m_n} = 2(1 - h_n^\tau)^{m_n} \leq 2 \exp(-m_n h_n^\tau).$$

Since we assume (4.32) and (H.iii), we see that, for all $\gamma > 1$ and $n$ large, $P_n \leq \exp(-(\log n)^\gamma)$, from which we readily conclude (4.30) by the Borel–Cantelli lemma. The case $I(\vartheta) = 0$ is readily inferred from the $0 < 2I(\vartheta) \leq 1$ case. $\square$

4.2. *Proof of part* (a) *of Theorem* 1.

LEMMA 3. *For some constant $C > 0$ independent of the sequence $\{h_n\}_{n \geq 1}$, with probability* 1,

(4.33) $$\limsup_{n \to \infty} \sup_{z \in J} \sup_{g \in \mathcal{G}} |L_n(z, g)| \leq C.$$

PROOF. The proof will be obtained by blocking between $2^k$ and $2^{k+1}$ and using Inequality 1. Notice that for a suitable $\tau > 0$, for all large $k$,

$$\max_{2^k < n \leq 2^{k+1}} \sup_{z \in J} \sup_{g \in \mathcal{G}} Eg^2(h_n^{-1/d}(z - Z)) \leq \tau^2 h_{n_{k+1}},$$

with $n_k := 2^k$, $k = 1, 2, \ldots$. Set for $k = 1, 2, \ldots$,

$$\mathcal{F}_{n_{k+1}} = \{g(h_n^{-1/d}(z - \cdot)) : g \in \mathcal{G},\ z \in J,\ 2^k < n \leq 2^{k+1}\}.$$

Now for any $\rho > 0$ and $D_0$ as in Inequality 1 with $\beta_1$ as in (4.18), we get, using $h_n \log(1/h_n) \searrow$, that

$$P\left\{\max_{2^k < n \leq 2^{k+1}} \sup_{z \in J} \sup_{g \in \mathcal{G}} |L_n(z, g)| > \sqrt{2}(\tau + \rho) D_1/\sqrt{\beta_1}\right\}$$

$$\leq P\left\{\max_{1 \leq n \leq n_{k+1}} \|T_n\|_{\mathcal{F}_{n_{k+1}}} \geq (\tau + \rho) D_1 b_{n_{k+1}}\right\},$$

which since $\mathcal{F}_{n_{k+1}}$ satisfies (F.ii–iii), is, by Inequality 1,

$$\leq 2 \exp(-D_0(\rho/\tau)^2 \log(1/h_{n_{k+1}})).$$

Notice that by (H.iii), we have $\log(1/h_{n_{k+1}})/\log(k) \to \infty$, which in combination with this last bound and the Borel–Cantelli lemma implies that (4.33) holds with $C = \sqrt{2}(\tau + \rho) D_1/\sqrt{\beta_1}$ for any $\rho > 0$. $\square$

Write now, for any $\gamma > 0$,

(4.34) $$\nu_k = [(1 + \gamma)^k] \qquad \text{for } k = 1, 2, \ldots.$$

LEMMA 4. *With probability 1,*

(4.35) $$\lim_{\gamma \searrow 0} \limsup_{k \to \infty} \max_{\nu_k < n \leq \nu_{k+1}} |b_n/b_{\nu_{k+1}} - 1| \sup_{z \in J} \sup_{g \in \mathcal{G}} |L_n(z, g)| = 0$$

*and*

(4.36) $$\lim_{\gamma \searrow 0} \limsup_{k \to \infty} \max_{\nu_k < n \leq \nu_{k+1}} \sup_{z \in J} \sup_{g \in \mathcal{G}} (b_n/b_{\nu_{k+1}}) |L_n(z, g((h_n/h_{\nu_{k+1}})^{1/d} \cdot))$$
$$- L_n(z, g)| = 0.$$

PROOF. The proofs of (4.35) and (4.36) follow closely those of Lemmas 3.5 and 3.6 of Deheuvels and Mason (1992). Lemma 3 is used to establish (4.35). The proof of (4.35) is based on Inequality 1 using condition (G.ii). We omit the routine details. $\square$

By condition (G.i) and compactness of $J$, for any $0 < \theta < 1$, we can choose $z_1, \ldots, z_{M_n(\theta)} \in J$ with $M_n(\theta) < \infty$ such that, for all $z \in J$,

$$\sup_{g \in \mathcal{G}} \min_{1 \leq i \leq M_n(\theta)} \int_{\mathbb{R}^d} [g(x) - g(x + h_n^{-1/d}(z_i - z))]^2 \, dx \leq \theta$$

and, further, we can do this so that

(4.37) $$\sup_{z \in J} \min_{1 \leq i \leq M_n(\theta)} h_n^{-1/d} |z - z_i|_2 \to 0 \quad \text{as } \theta \searrow 0,$$

and for some function $C(\theta) < \infty$ for $\theta > 0$,

(4.38) $$M_n(\theta) \leq C(\theta) h_n^{-1}.$$

Next, for any $0 < \theta < 1$, $z \in J$, let $z(\theta)$ denote a selection of a $z_i$ among $z_1, \ldots, z_{M_n(\theta)}$ satisfying

(4.39) $$\sup_{g \in \mathcal{G}} \int_{\mathbb{R}^d} [g(x) - g(x + h_n^{-1/d}(z(\theta) - z))]^2 \, dx \leq \theta.$$

Moreover, we do this in such a way so that

(4.40) $$\sup_{z \in J} h_n^{-1/d} |z - z(\theta)|_2 \to 0 \quad \text{as } \theta \searrow 0.$$

LEMMA 5. *There exists a $\tau > 0$ such that, for all $0 < \theta < 1$ and large enough $n$,*

(4.41) $$\sup_{z \in J} \sup_{g \in \mathcal{G}} E(g(h_n^{-1/d}(z - Z)) - g(h_n^{-1/d}(z(\theta) - Z)))^2 \leq \tau^2 \theta h_n.$$

PROOF. Notice that

$$E(g(h_n^{-1/d}(z - Z)) - g(h_n^{-1/d}(z(\theta) - Z)))^2$$
$$\leq E[E(g(h_n^{-1/d}(z - Z)) - g(h_n^{-1/d}(z(\theta) - Z)))^2 | Z \in \Omega_n(z, z(\theta))],$$

where $\Omega_n(z, z(\theta)) = (z - h_n^{1/d} I^d) \cup (z(\theta) - h_n^{1/d} I^d)$. Now this last bound equals

$$\frac{\int_{\Omega_n(z,z(\theta))} [g(h_n^{-1/d}(z - y)) - g(h_n^{-1/d}(z(\theta) - y))]^2 f(y) \, dy}{P(\Omega_n(z, z(\theta)))} P(\Omega_n(z, z(\theta)))$$
$$=: I_n(z, z(\theta)).$$

Clearly, for all large enough $n$, uniformly in $z \in J$, $\Omega_n(z, z(\theta)) \subset J_\gamma$. Thus by using the fact that, with $\lambda$ denoting Lebesgue measure, $h_n \leq \lambda(\Omega_n(z, z(\theta))) \leq$

$2h_2$, along with (4.18), we get, for all large enough $n$, uniformly in $z \in J$ and $g \in \mathcal{G}$,

$$I_n(z, z(\theta)) \leq \frac{\int_{\Omega_n(z,z(\theta))}[g(h_n^{-1/d}(z-y)) - g(h_n^{-1/d}(z(\theta)-y))]^2 \beta_2 \, dy}{h_n \beta_1} 2\beta_2 h_n$$

$$\leq \frac{\int_{\mathbb{R}^d}[g(h_n^{-1/d}(z-y)) - g(h_n^{-1/d}(z(\theta)-y))]^2 \beta_2 \, dy}{h_n \beta_1} 2\beta_2 h_n,$$

which by the change of variables $x = h_n^{-1/d}(z-y)$ and (4.39)

$$= \frac{\int_{\mathbb{R}^d}[g(x) - g(x + h_n^{-1/d}(z(\theta)-z))]^2 \, dx}{\beta_1} 2\beta_2^2 h_n \leq 2\beta_2^2 \theta h_n / \beta_1.$$

Thus we have (4.41) with $\tau^2 = 2\beta_2^2/\beta_1$.  $\square$

For each $z \in J$ and $h > 0$, let $S_n(\,\cdot\,; z, h)$ denote the function from $\mathcal{G}$ to $\mathbb{R}$ defined, for each $g \in \mathcal{G}$, to be

$$(4.42) \quad S_n(g; z, h) = (\sqrt{f(z)})^{-1} \sum_{j=1}^{n} \{g(h^{-1/d}(z - Z_j)) - Eg(h^{-1/d}(z - Z))\}.$$

Now for any $0 < \theta < 1$, with $z(\theta)$ and $M_n(\theta)$ satisfying (4.37)–(4.40), set

$$\varpi_n^{(1)}(\theta) = b_n^{-1} \sup_{z \in J} \sup_{g \in \mathcal{G}} |S_n(g; z, h_n) - \sqrt{f(z(\theta))/f(z)} S_n(g; z(\theta), h_n)|$$

and

$$\varpi_n^{(2)}(\theta) = \delta(\theta) b_n^{-1} \sup_{z \in J} \sup_{g \in \mathcal{G}} |S_n(g; z(\theta), h_n)|,$$

where $b_n$ is as in (4.1) and $\delta(\theta) = \sup_{z \in J}|\sqrt{f(z(\theta))/f(z)} - 1|$. Notice that

$$\sup_{z \in J} \sup_{g \in \mathcal{G}} |L_n(z, g) - L_n(z(\theta), g)| \leq \varpi_n^{(1)}(\theta) + \varpi_n^{(2)}(\theta),$$

from which we get that, for any $\eta > 0$ and $0 < \theta < 1$,

$$P\{L_n(z, \cdot) \notin \mathcal{S}_0^\eta \text{ for some } z \in J\}$$

$$(4.43) \quad \leq \sum_{i=1}^{M_n(\theta)} P\{L_n(z_i, \cdot) \notin \mathcal{S}_0^{\eta/2}\}$$

$$+ P\{\varpi_n^{(1)}(\theta) > \eta/4\} + P\{\varpi_n^{(2)}(\theta) > \eta/4\}.$$

To complete the proof of part (a) of Theorem 1, we need the following generalization of the Ottaviani inequality.

4.2.1. *A generalized Ottaviani inequality.* Let $\{S_m(t):t\in\Lambda, 0\leq m\leq n\}$, $n\geq 1$, be an indexed set of random processes such that, for each $t\in\Lambda$ and $1\leq m\leq n$, $S_m(t)\in B$, where $\Lambda$ is a countable set and $B$ is a separable Banach space with norm $\|\cdot\|$. Also assume that, for each $1\leq m\leq n$,

(4.44) $\{S_n(t) - S_m(t):t\in\Lambda\}$ is independent of $\{S_k(t):t\in\Lambda, 1\leq k\leq m\}$.

Further assume that, for some $\tau > 0$,

(4.45) $$\max_{0\leq m\leq n}\sup_{t\in\Lambda} P\{\|S_n(t) - S_m(t)\| \geq \tau\} =: c < 1,$$

where $S_0(t) = 0$ for all $t\in\Lambda$. For any Borel subset $A\subset B$, $n\geq 1$ and $\delta > 0$, set $A^\delta = \{x:\inf_{y\in A}\|x-y\| < \delta\}$,

$$C_n(\delta) = \{S_m(t)\notin A^\delta \text{ for some } t\in\Lambda \text{ and } 1\leq m\leq n\}$$

and

$$D_n(\delta) = \{S_n(t)\notin A^\delta \text{ for some } t\in\Lambda\}.$$

We shall prove the following extension of Ottaviani's inequality.

LEMMA 6. *With $\tau > 0$ and $0\leq c < 1$ as in* (4.45), *for all Borel subsets $A\subset B$, and $\varepsilon > 0$,*

(4.46) $$P\{C_n(\varepsilon + \tau)\} \leq (1-c)^{-1} P\{D_n(\varepsilon)\}.$$

PROOF. Let $t_i$, $i=1,2,\ldots$, be an enumeration of the set $\Lambda$ and define the events for $\delta > 0$, $i=1,2,\ldots$, and $1\leq m\leq n$,

$$D_{i,m}(\delta) = \{S_m(t_i)\notin A^\delta\}, \qquad D_m(\delta) = \bigcup_{i\geq 1} D_{i,m}(\delta)$$

and

$$F_{i,m}(\delta) = \{\|S_n(t_i) - S_m(t_i)\| < \delta\}.$$

We define $D_{i,0}(\delta) = \varnothing$ and $D_0(\delta) = \varnothing$. We get that

$$P\{C_n(\varepsilon+\tau)\} = \sum_{q=1}^n \sum_{i=1}^\infty P\bigg\{D_{i,q}(\varepsilon+\tau) \bigcap_{j\leq i-1} D_{j,q}^C(\varepsilon+\tau) \bigcap_{k\leq q-1} D_k^C(\varepsilon+\tau)\bigg\},$$

where $A^C$ denotes the complement of the event $A$. We see by (4.44) that, for any $\delta'$ and $\delta'' > 0$, the sequences of events $\{D_{i,m}(\delta'):i\geq 1\}$ and $\{F_{i,m}(\delta''):i\geq 1\}$ are independent. Therefore by (4.45),

$(1-c)P\{C_n(\varepsilon+\tau)\}$
$\quad = \min_{1\leq m\leq n}\inf_{i\geq 1} P\{F_{i,m}(\tau)\}P\{C_n(\varepsilon+\tau)\}$

$$\leq \sum_{i=1}^{\infty} \sum_{q=1}^{n} P\left\{ D_{i,q}(\varepsilon+\tau) \cap F_{i,q}(\tau) \bigcap_{j \leq i-1} D_{j,q}^C(\varepsilon+\tau) \bigcap_{k \leq q-1} D_k^C(\varepsilon+\tau) \right\}$$

$$\leq P\left\{ \bigcup_{i=1}^{\infty} D_{i,n}(\varepsilon) \right\} = P\{D_n(\varepsilon)\}. \qquad \square$$

The ideas used in the proof of Lemma 6 go back at least to Lemma 2.3 of James (1975).

Returning now to the proof of part (a) of Theorem 1 and recalling (4.34), consider, for any $\varepsilon, \gamma > 0$ and $k \geq 1$, the sets

$$C_k(\varepsilon, \gamma) = \{b_{\nu_{k+1}}^{-1} S_n(\,\cdot\,; z, h_{\nu_{k+1}}) \notin \mathcal{S}_0^{\varepsilon} \text{ for some } \nu_k < n \leq \nu_{k+1} \text{ and } z \in J\}$$

and

$$D_k(\varepsilon, \gamma) = \{b_{\nu_{k+1}}^{-1} S_{\nu_{k+1}}(\,\cdot\,; z, h_{\nu_{k+1}}) \notin \mathcal{S}_0^{\varepsilon} \text{ for some } z \in J\}.$$

It is elementary to verify using Inequality 1 that, for any $\varepsilon > 0$, all $\gamma > 0$ small, depending on $\varepsilon$, and all large enough $k$,

$$\max_{\nu_k < m \leq \nu_{k+1}} \sup_{z \in J} P\{b_{\nu_{k+1}}^{-1} \|S_{\nu_{k+1}}(\,\cdot\,; z, h_{\nu_{k+1}}) - S_m(\,\cdot\,; z, h_{\nu_{k+1}})\|_{\mathcal{F}} \geq \varepsilon/2\} < 1/2.$$

Thus, since with probability 1, the values of $S_n(\,\cdot\,; z, h_{\nu_{k+1}})$, $\nu_k < n \leq \nu_{k+1}$, $z \in J$, are determined by a countable subset of $J$, it is clear that we can apply the generalized Ottaviani inequality to give, for all large $k$,

$$P\{C_k(\varepsilon, \gamma)\} \leq 2P\{D_k(\varepsilon/2, \gamma)\} = 2P\{L_{\nu_{k+1}}(z, \cdot) \notin \mathcal{S}_0^{\varepsilon/2} \text{ for some } z \in J\},$$

which by inequality (4.43) is

$$\leq 2 \sum_{i=1}^{M_{\nu_{k+1}}(\theta)} P\{L_{\nu_{k+1}}(z_i, \cdot) \notin \mathcal{S}_0^{\varepsilon/4}\} + 2P\{\varpi_{\nu_{k+1}}^{(1)}(\theta) > \varepsilon/8\}$$

$$+ 2P\{\varpi_{\nu_{k+1}}^{(2)}(\theta) > \varepsilon/8\}$$

$$=: Q_{1,k}(\varepsilon) + Q_{2,k}(\varepsilon) + Q_{3,k}(\varepsilon).$$

First, by choosing $\theta > 0$ sufficiently small and using the fact that (4.40) implies that $\delta(\theta) \to 0$, as $\theta \searrow 0$, along with (4.41), one can easily show using Inequality 1 that

$$\sum_{k=1}^{\infty} (Q_{2,k}(\varepsilon) + Q_{3,k}(\varepsilon)) < \infty.$$

Next, by applying part (i) of Proposition 1 with $m_n = M_{\nu_{k+1}}(\theta)$, in combination with Fact 6 with $m_n = 1$, it is straightforward to check that, for some $\eta > 0$ and all large $k$,

$$Q_{1,k}(\varepsilon) \leq M_{\nu_{k+1}}(\theta) \exp(-(1+\eta) \log(1/h_{\nu_{k+1}})),$$

which by (4.38) and $\log(1/h_{\nu_{k+1}})/\log(k) \to \infty$, following from (H.iii), implies that $\sum_{k=1}^{\infty} Q_{1,k}(\varepsilon) < \infty$. Thus for any $\varepsilon > 0$ and all $\gamma > 0$ small, depending on $\varepsilon$, using the Borel–Cantelli lemma and the above string of inequalities, we obtain that

(4.47) $$P\{C_k(\varepsilon, \gamma), \text{ i.o. in } k \geq 1\} = 0.$$

Observing that for $\nu_k < n \leq \nu_{k+1}$,

$$(b_n/b_{\nu_{k+1}})L_n(z,\ g((h_n/h_{\nu_{k+1}})^{1/d} \cdot)) = S_n(\,\cdot\,; z, h_{\nu_{k+1}})/b_{\nu_{k+1}},$$

the remainder of the proof of part (a) is now easily inferred from (4.47) and Lemma 4.

4.3. *Proof of Corollary* 1. First note that the assumption that the density $f$ is uniformly continuous on $\mathbb{R}^d$ is equivalent to the assumption that $f$ is continuous on $\mathbb{R}^d$ and satisfies the condition that

$$\lim_{R \to \infty} \sup\{f(z) : |z|_2 \geq R\} = 0,$$

from which we readily infer that $\tau_0 = \sup_{z \in \mathbb{R}^d} \sqrt{f(z)} < \infty$. Furthermore, $f$ continuous on $\mathbb{R}^d$ implies that, for all $c > 0$, the set $\{z : c > f(z) > 0\}$ is nonempty.

Define the compact set $J = \{z : c \leq f(z) \text{ and } |z|_2 \leq 2R\}$, where $c > 0$ and $R$ are chosen so that $J$ has nonempty interior and with $D_1$ as in Inequality 1 and $\kappa$ the bound on the functions in $\mathcal{G}$,

(4.48) $$\sup\{\sqrt{f(z)} : |z|_2 \geq R\} < \sqrt{c} \leq \tau_0/(6\sqrt{2}\kappa D_1).$$

Now since $f$ is assumed to be uniformly continuous, we can choose a $\gamma > 0$ so that

$$f(z) \geq c/2 \quad \text{for all } z \in J_\gamma.$$

Thus we can apply Theorem 1 to conclude that, for all $\varepsilon > 0$, there exists an $n(\varepsilon)$ such that for each $n \geq n(\varepsilon)$, $\{L_n(z, \cdot) : z \in J\} \subset \mathcal{S}_0^\varepsilon$, which clearly implies that $\{D_n(z, \cdot) : z \in J\} \subset \tau_0 \mathcal{S}_0^\varepsilon$. Obviously now, to complete the proof of the first part of Corollary 1, it suffices to show that

(4.49) $$\limsup_{n \to \infty} \sup_{z \in B_R} \sup_{g \in \mathcal{G}} |D_n(z, g)| \leq \tau_0/2 \quad \text{a.s.,}$$

where $B_R = \{z : |z|_2 \geq R\}$. The proof will follow from blocking between $2^k$ and $2^{k+1}$ and using Inequality 1. Notice that since each $g \in \mathcal{G}$ is bounded by $\kappa > 0$, we get, for each $n \geq 1$, $z \in B_R$ and $g \in \mathcal{G}$,

$$Eg^2(h_n^{-1/d}(z - Z)) \leq \kappa^2 h_n f * H_{h_n}(z),$$

FUNCTIONAL LAW 27

where $H(x) = I\{x \in I^d\}$, which by Lemma 1, for all $n$ large enough uniformly in $z \in B_R$, is

$$\leq \kappa^2 h_n(f(z) + c) \leq 2c\kappa^2 h_n.$$

Thus with $n_k := 2^k$, $k = 1, 2, \ldots$, using (H.i), we get, for all large enough $k \geq 1$,

$$\max_{2^k < n \leq 2^{k+1}} \sup_{z \in B_R} \sup_{g \in \mathcal{G}} E g^2(h_n^{-1/d}(z - Z)) \leq 4c\kappa^2 h_{n_{k+1}}.$$

Set, for $k = 1, 2, \ldots,$

$$\mathcal{F}_{n_{k+1}} = \{g(h_n^{-1/d}(z - \cdot)) : g \in \mathcal{G}, \ z \in B_R, \ 2^k < n \leq 2^{k+1}\}.$$

Now with $D_1$ as in Inequality 1,

$$P\left\{\max_{2^k < n \leq 2^{k+1}} \sup_{z \in B_R} \sup_{g \in \mathcal{G}} |D_n(z, g)| > \sqrt{2}(2\kappa\sqrt{c} + \kappa\sqrt{c})D_1\right\}$$

$$\leq P\left\{\max_{1 \leq n \leq n_{k+1}} \|T_n\|_{\mathcal{F}_{n_{k+1}}} \geq (2\kappa\sqrt{c} + \kappa\sqrt{c})D_1 b_{n_{k+1}}\right\},$$

which since $\mathcal{F}_{n_{k+1}}$ satisfies (F.ii–iii) is, by Inequality 1, with $\rho = \kappa\sqrt{c}$ and $\gamma = 2\kappa\sqrt{c}$,

$$\leq 2\exp(-D_0 4^{-1} \log(1/h_{n_{k+1}})).$$

Since $\log(1/h_{n_{k+1}})/\log(k) \to \infty$, this last bound, the Borel–Cantelli lemma and (4.48) imply that (4.49) holds, which completes the proof of part (a) of Corollary 1.

Whenever $f(z) > 0$, part (b) of Corollary 1 is proved by applying part (b) of Theorem 1 on closed balls $\overline{B}_\delta(z) = \{x : |x - z|_2 \leq \delta\}$ of radius $\delta > 0$ around $z$, where $\delta > 0$ is sufficiently small, and when $f(z) = 0$ we apply a straightforward modification of the argument given in the previous paragraph to closed balls $\overline{B}_\delta(z)$ around $z$ to show that

$$\lim_{\delta \searrow 0} \limsup_{n \to \infty} \sup_{z \in \overline{B}_\delta(z)} \sup_{g \in \mathcal{G}} |D_n(z, g)| = 0 \quad \text{a.s.}$$

**Acknowledgments.** The author thanks Alfredo Arreguin, Paul Eggermont, Uwe Einmahl, Evarist Giné, David Pollard and Andre Zaitsev for useful discussions while this paper was being written. In addition, he is especially grateful to Philippe Berthet and Davit Varron for their careful reading of the manuscript and to Miguel Arcones for sending him a preliminary version of his work on large deviations.

## REFERENCES

Alexander, K. S. (1987). Rates of growth and sample moduli for weighted empirical processes. *Probab. Theory Related Fields* **75** 379–423. MR890285

Arcones, M. (2003). The large deviation principle of stochastic processes. I. *Theory Probab. Appl.* **47** 567–583. MR2001788

Arcones, M. (2004). The large deviation principle of stochastic processes. II. *Theory Probab. Appl.* **48** 19–44.

Borisov, I. (2002). Moment inequalities connected with accompanying Poisson laws in Abelian groups. Preprint. MR2003788

Deheuvels, P. (2000). Uniform limit laws for kernel density estimators on possibly unbounded intervals. In *Recent Advances in Reliability Theory: Methodology, Practice and Inference* (N. Limnios and M. Nikulin, eds.) 477–492. Birkhäuser, Basel. MR1783500

Deheuvels, P. and Mason, D. M. (1992). Functional laws of the iterated logarithm for increments of empirical and quantile processes. *Ann. Probab.* **20** 1248–1287. MR1175262

Deheuvels, P. and Mason, D. M. (1994). Functional laws of the iterated logarithm for local empirical processes indexed by sets. *Ann. Probab.* **22** 1619–1661. MR1303659

Deheuvels, P. and Mason, D. M. (2004). General asymptotic confidence bounds based on kernel–type function estimators. *Stat. Inference Stoch. Process.* To appear.

Devroye, L. and Lugosi, G. (2000). *Combinatorial Methods in Density Estimation.* Springer, New York. MR1843146

Einmahl, J. H. J. (1987). *Multivariate Empirical Processes. CWI Tract* **32**. Centrum Wisk. Inform., Amsterdam. MR890957

Einmahl, U. and Mason, D. M. (1997). Gaussian approximation of local empirical processes indexed by functions. *Probab. Theory Related Fields* **107** 283–311. MR1440134

Einmahl, U. and Mason, D. M. (1998). Strong approximations for local empirical processes. In *Proceedings of High Dimensional Probability* (E. Eberlein, M. Hahn and J. Kuelbs, eds.) 75–92. Birkhäuser, Basel. MR1652321

Einmahl, U. and Mason, D. M. (2000). An empirical process approach to the uniform consistency of kernel-type function estimators. *J. Theoret. Probab.* **13** 1–37. MR1744994

Einmahl, U. and Mason, D. M. (2003). Uniform in bandwidth consistency of kernel-type function estimators. Preprint.

Giné, E. and Guillou, A. (2001). On consistency of kernel density estimators for randomly censored data: Rates holding uniformly over adaptive intervals. *Ann. Inst. H. Poincaré Probab. Statist.* **37** 503–522. MR1876841

Giné, E. and Guillou, A. (2002). Rates of strong consistency for multivariate kernel density estimators. *Ann. Inst. H. Poincaré Probab. Statist.* **38** 907–921. MR1955344

Giné, E., Koltchinskii, V. and J. Wellner, J. A. (2003). Ratio limit theorems for empirical processes. In *Stochastic Inequalities and Applications* 249–278. Birkhäuser, Basel. MR2073436

Giné, E., Mason, D. M. and Zaitsev, A. (2003). The $L_1$-norm density estimator process. *Ann. Probab.* **31** 719–768. MR1964947

Giné, E. and Zinn, J. (1984). Some limit theorems for empirical processes. *Ann. Probab.* **12** 929–989. MR757767

Härdle, W., Janssen, P. and Serfling, R. (1988). Strong uniform consistency rates for estimators of conditional functionals. *Ann. Statist.* **16** 1428–1449. MR964932

James, B. (1975). A functional law of the iterated logarithm for weighted empirical distributions. *Ann. Probab.* **3** 762–772. MR402881

Komlós, J., Major, P. and Tusnády, G. (1975). An approximation of partial sums of independent rv's and the sample df I. *Z. Wahrsch. Verw. Gebiete* **32** 111–131. MR375412

Mason, D. M. (1988). A strong invariance theorem for the tail empirical process. *Ann. Inst. H. Poincaré Sect. B* **24** 491–506. MR978022

Mason, D. M. (2003). A uniform functional law of the logarithm for a local Gaussian process. In *High Dimensional Probability III* (J. Hoffmann-Jorgensen, M. B. Marcus and J. A. Wellner, eds.) 135–151. Birkhäuser, Boston. MR2033886

Nolan, D. and Marron, J. S. (1989). Uniform consistency of automatic and location-adaptive delta-sequence estimators. *Probab. Theory Related Fields* **80** 619–632. MR980690

Nolan, D. and Pollard, D. (1987). $U$-processes: Rates of convergence. *Ann. Statist.* **15** 780–799. MR888439

Parzen, E. (1961). An approach to time series analysis. *Ann. Math. Statist.* **32** 951–989. MR143315

Pyke, R. and Shorack, G. R. (1968). Weak convergence of a two–sample empirical process and a new approach to Chernoff–Savage theorems. *Ann. Math. Statist.* **39** 755–771. MR226770

Rio, E. (1994). Local invariance principles and their applications to density estimation. *Probab. Theory Related Fields* **98** 21–45. MR1254823

Stein, E. M. (1970). *Singular Integrals and Differentiability Properties of Functions*. Princeton Univ. Press. MR290095

Stoyan, D., Kendall, W. S. and Mecke, J. (1995). *Stochastic Geometry and Its Applications*. Wiley, New York. MR895588

Stute, W. (1982a). The oscillation behavior of empirical processes. *Ann. Probab.* **10** 86–107. MR637378

Stute, W. (1982b). The law of the iterated logarithm for kernel density estimators. *Ann. Probab.* **10** 414–422. MR647513

Stute, W. (1984). The oscillation behavior of empirical processes: The multivariate case. *Ann. Probab.* **12** 361–379. MR735843

Talagrand, M. (1994). Sharper bounds for Gaussian and empirical processes. *Ann. Probab.* **22** 28–76. MR1258865

van der Vaart, A. W. (1998). *Asymptotic Statistics*. Cambridge Univ. Press. MR1652247

van der Vaart, A. W. and Wellner, J. A. (1996). *Weak Convergence and Empirical Processes*. Springer, New York. MR1385671

Zaitsev, A. Yu. (1987a). Estimates of the Lévy–Prokhorov distance in the multivariate central limit theorem for random variables with finite exponential moments. *Theory Probab. Appl.* **31** 203–220.

Zaitsev, A. Yu. (1987b). On the Gaussian approximation of convolutions under multidimensional analogues of S. N. Bernstein's inequality conditions. *Probab. Theory Related Fields* **74** 534–566. MR876255

Statistics Program  
University of Delaware  
206 Townsend Hall  
Newark, Delaware 19716  
USA  
e-mail: davidm@udel.edu